\documentclass[12pt]{article}
\usepackage{amssymb,amsmath,mathrsfs}
\usepackage{hyperref}
\usepackage{graphicx}
\usepackage{color}
\usepackage[OT2,OT1]{fontenc}
\usepackage{upquote}
\usepackage{authblk}
\usepackage{mathtools}

\DeclarePairedDelimiter{\ceil}{\lceil}{\rceil}

\begin{document}

\title{\LARGE\bf A single-domain implementation \\
of the Voigt/complex error function \\
by vectorized interpolation}

\bigskip
\author[1, 2]{\small S. M. Abrarov}
\author[1, 3]{\small B. M. Quine}
\author[1, 2, 3]{\small R. Siddiqui}
\author[2, 3]{\small R. K. Jagpal}

\affil[1]{\scriptsize Dept. Earth and Space Science and Engineering, York University, 4700 Keele St., Canada, M3J 1P3 \normalsize}
\affil[2]{\scriptsize Epic College of Technology, 5670 McAdam Rd., Mississauga, Canada, L4Z 1T2 \normalsize}
\affil[3]{\scriptsize Dept. Physics and Astronomy, York University, 4700 Keele St., Toronto, Canada, M3J 1P3 \normalsize}

\date{August 31, 2019}
\maketitle

\begin{abstract}
In this work we show how to perform a rapid computation of the Voigt/complex error over a single domain by vectorized interpolation. This approach enables us to cover the entire set of the parameters $x,y \in \mathbb{R}$ required for the HITRAN-based spectroscopic applications. The computational test reveals that within domains $x\in\left[0,15\right]\cap y\in\left[10^{-8},15\right]$ and $x\in\left[0,50000\right]\cap y\geq 10^{-8}$ our algorithmic implementation is faster in computation by factors of about $8$ and $3$, respectively, as compared to the fastest known C/C++ code for the Voigt/complex error function. A rapid MATLAB code is presented.
\vspace{0.25cm}
\vspace{0.25cm}
\end{abstract}

\section{Introduction}
The complex error function, also commonly known as the Faddeeva function, can be defined as  \cite{Faddeeva1961, Armstrong1967, Gautschi1970, Abramowitz1972}
\begin{equation}\label{eq_1}
w\left( z \right) = {e^{ - {z^2}}}\left( {1 + \frac{{2i}}{{\sqrt \pi  }}\int\limits_0^z {{e^{{t^2}}}dt} } \right),
\end{equation}
where $z = x + iy$ is the complex argument. The real part of the complex error function \eqref{eq_1}, known as the Voigt function, can be represented as given by \cite{Armstrong1967, Srivastava1987, Srivastava1992}
\[
\label{eq_2a}\tag{2a}
K\left( {x,y} \right) = \frac{1}{{\sqrt \pi  }}\int\limits_0^\infty  {{e^{ - {t^2/4}}}{e^{ - yt}}\cos \left( {xt} \right)dt}
\]
that is proportional to the Voigt profile describing the spectral broadening of atmospheric gas absorption or emission
$$
g_V\left(\nu-\nu_0,\alpha_L,\alpha_D\right)=\frac{\sqrt{\ln 2/\pi}}{\alpha_D}K\left(x,y\right),
$$
where $\nu$ is the frequency, $\nu_0$ is the frequency of the line center, $\alpha_L$ and $\alpha_D$ are the Lorentz and Doppler half widths at half maximum (HWHM), respectively, and
$$
x=\sqrt{\ln 2}\,\frac{\nu-\nu_0}{\alpha_D}, \qquad\qquad y=\sqrt{\ln 2}\,\frac{\alpha_L}{\alpha_D}.
$$

The imaginary part of the complex error function \eqref{eq_1} has no specific name. Historically, it is denoted by $L\left( {x,y} \right)$ and can be written in form \cite{Srivastava1987, Srivastava1992}
\[
\label{eq_2b}\tag{2b}
L\left( {x,y} \right) = \frac{1}{{\sqrt \pi  }}\int\limits_0^\infty  {{e^{ - {t^2/4}}}{e^{ - yt}}\sin \left( {xt} \right)dt}.
\]

None of the integrals above are analytically integrable in closed form. Therefore, these equations must be solved numerically. There are two most important aspects that have to be taken into consideration for efficient algorithmic implementation of the integrals above in spectroscopic applications based on HITRAN database \cite{Hill2016}. Specifically, the latest versions of the HITRAN database provide spectroscopic data with $4$ and more digits in floating format. In order to exclude the truncation errors that inevitably occur in any line-by-line atmospheric modeling \cite{Berk2017, Pliutau2017, Siddiqui2015, Siddiqui2017, Jagpal2010, Quine2013}, the algorithmic implementation of the Voigt/complex error function should provide accuracy at least as close as possible to ${10^{ - 6}}$. This accuracy requirement is particularly relevant for the set of input numbers $\left\{x,y \in \mathbb{R}\hspace{-0.1cm}:\,\left|x\right| + y\leq 15 \right\}$ while it is not so much critical for the set $\left\{x,y \in \mathbb{R}\hspace{-0.1cm}:\,\left|x\right| + y > 15 \right\}$ according to literature \cite{Schreier2011, Schreier2018}. Furthermore, the algorithm should be rapid as it may deal with many millions of input numbers $z$ in computation for the various radiative transfer applications \cite{Schreier2018, Grimm2015}.

In order to satisfy these two criteria in computation of the integrals above, the complex plane is segmented into several domains. For example, Huml\'{i}\v{c}ek proposed a rapid algorithm for computation of the Voigt function \eqref{eq_2a} based on rational functions by segmenting the complex plane into several domains such that each of them is computed by corresponding rational approximation \cite{Humlicek1982}. Kuntz proposed some modifications and showed efficiency of the Huml\'{i}\v{c}ek algorithm with $4$ domains \cite{Kuntz1997}. The Huml\'{i}\v{c}ek's algorithm is rapid and provides accuracy ${10^{ - 4}}$. Later, Wells developed a FORTRAN code where he succeeded to improve accuracy by an order of the magnitude by making some modifications to the original Huml\'{i}\v{c}ek's algorithm \cite{Wells1999}. Another interesting variation of the Huml\'{i}\v{c}ek's algorithm with improved accuracy was presented by Imai {\it{et al.}} \cite{Imai2010}.

It is very desirable to reduce wherever possible the number of the domains to accelerate an algorithm since each area segmentation of the complex plane leads to run-time increase due to additional logical operations and subsequent sorting of the input numbers $x$ and $y$.

In our earlier publication \cite{Abrarov2009} we have shown how to reduce the number of domains to $2$ by interpolation. In particular, we suggested two-domain scheme for rapid computation of the Voigt function that can be represented as
\setcounter{equation}{2}
\small
\begin{equation}\label{eq_3}
K\left( {x,y} \right) \approx \left\{ 
\begin{aligned}
&{\text{interpolation,}} \qquad \hspace{0.43cm}\frac{x^2}{27^2} + \frac{y^2}{15^2} \leqslant 1 \\
&\frac{{{a_1} + {b_1}{x^2}}}{{{a_2} + {b_2}{x^2} + {x^4}}}, \qquad \frac{x^2}{27^2} + \frac{y^2}{15^2} > 1,
\end{aligned}\right.
\end{equation}
\normalsize
where (see Appendix in \cite{Kuntz1997})
\small
\[
\begin{aligned}
&a_1 = y/(2\sqrt{\pi})+y^3/\sqrt{\pi}\approx 0.2820948y + 0.5641896y^3 \\
&b_1 = y/\sqrt{\pi}\approx 0.5641896y \\
&a_2 = 0.25 + y^2 + y^4 \\
&b_2 =-1 + 2y^2
\end{aligned}
\]
\normalsize
and
$$
\frac{{{a_1} + {b_1}{x^2}}}{{{a_2} + {b_2}{x^2} + {x^4}}}=\operatorname{Re} \left\{ \frac{iz/\sqrt{\pi}}{z^2-1/2}\right\}.
$$
Thus, the complex plane is segmented into two domains by an ellipse centered at the origin with semi-major and semi-minor axises equal to $27$ and $15$, respectively (see inset in Fig. 2 from our paper \cite{Abrarov2009}). Inside the ellipse (computationally difficult internal domain) we apply interpolation while outside the ellipse (computationally simple external domain) we apply a simple rational approximation of low order. This scheme shows high efficiency particularly when ${\mathbf{x}} = \left\{ {{x_1},{x_2},{x_3}, \ldots } \right\}$ is a vector and $y$ is a scalar. In fact, vectorized ${\mathbf{x}} = \left\{ {{x_1},{x_2},{x_3}, \ldots } \right\}$ and scalar $y$ is a quite common technique in radiative transfer applications \cite{Lynas-Gray1993, Letchworth2007, Schreier2008}.

Recently Schreier reported a two-domain scheme (see equation (12) in \cite{Schreier2018})
\small
\begin{equation}\label{eq_4}
w\left(z\right) =K\left(x,y\right)+iL\left(x,y\right)\approx \left\{
\begin{aligned}
&\frac{\sum_{k=0}^{n-1}{\alpha_k z^k}}{\sum_{\ell=0}^{n}{\beta_\ell z^\ell}}, \qquad\qquad\left|x\right| + y \leq 15 \\
&\frac{iz/\sqrt{\pi}}{z^2-1/2}, \qquad\qquad \hspace{0.4cm}\left|x\right| + y > 15,
\end{aligned}\right.
\end{equation}
\normalsize
where $n$ is assumed to be an even integer, $\alpha_k$ and $\beta_\ell$ are the expansion coefficients that can be readily generated by Computer Algebra System (CAS) supporting symbolic programming. It is suggested that $n=20$ in the single quotient is sufficient for line-by-line calculations in atmospheric modeling \cite{Schreier2018}.

In general, for arbitrary $n$ (even or odd) representation of the Huml\'{i}\v{c}ek's approximation as a single quotient can be made by using notation $\ceil{\ldots}$ for the ceiling function as
\begin{equation}\label{eq_5}
w\left(z\right)\approx \frac{1}{2}\sum_{k=1}^{n}\left(\frac{\gamma_k+i\theta_k}{z-x_k+i\delta}-\frac{\gamma_k-i\theta_k}{z+x_k+i\delta}\right) = \frac{\sum_{k=0}^{2\ceil{n/2}-1}{\alpha_k z^k}}{\sum_{\ell=0}^{2\ceil{n/2}}{\beta_\ell z^\ell}},
\end{equation}
where
\[
\begin{aligned}
\gamma_k=&-\frac{1}{\pi}\omega_k e^{\delta^2}\sin\left(2x_k \delta\right), \\
\theta_k=&\frac{1}{\pi}\omega_k e^{\delta^2}\cos\left(2x_k \delta\right),
\end{aligned}
\]
$x_k$ are roots of the Hermite polynomial $H_n\left(x\right)$ of degree $n$ that can be defined by recurrence relations \cite{Weisstein1}
\[
H_{n+1}=2xH_n\left(x\right)-2n H_{n-1}\left(x\right), \qquad H_0\left(x\right)=1, \qquad H_1\left(x\right)=2x,
\]
$\delta$ is a fitting parameter that at $n=20$ can be taken as $1.55$ \cite{Schreier2018} and
$$
\omega_k = \frac{2^{n-1}n!\sqrt{\pi}}{n^2 H_{n-1}^2\left(x_k\right)}
$$
are weights of the Hermite polynomial $H_n\left(x\right)$ of degree $n$ \cite{Weisstein2}.

It is interesting to note that if integer $n$ is even and roots are given in ascending order such that $x_{k-1} < x_{k}$, then number of the summation terms in the Huml\'{i}\v{c}ek's approximation can be reduced by a factor of two (compare equations (27a) and (26b) from \cite{Berk2017})
\[
w\left(z\right)\approx \sum_{k=1}^{n/2}\left(\frac{\gamma_k+i\theta_k}{z-x_k+i\delta}-\frac{\gamma_k-i\theta_k}{z+x_k+i\delta}\right).
\]

The single quotient reformulation shown in equations \eqref{eq_4} and \eqref{eq_5} is interesting. However, by performing computational test we found empirically that deterioration of accuracy with decreasing $y$ is especially inherent to the single quotient reformulation of the Huml\'{i}\v{c}ek's approximation \eqref{eq_5} (deterioration of accuracy with decreasing $y$ is a common problem in computation of the Voigt/complex error function \cite{Armstrong1967, Amamou2013, Abrarov2011}). Although a multiple precision arithmetic may be used to resolve this problem, it needs a special package (see for example \cite{Tsarapkina2014}) that may affect computational speed and makes the MATLAB code inconvenient in practical applications.

In plasma physics of rarefied gases or at low atmospheric pressure that takes place in stratosphere, mesosphere and thermoshpere of the Earth, where the Doppler broadening considerably predominates over the Lorentz broadening, the value of $y$ dependent on the pressure and temperature may be relatively close to zero. This is particularly important as the latest versions of the HITRAN supply parameters for high temperatures almost reaching $10000$ K and there is a tendency that it will be increased in future. Therefore, it would be very desirable to develop a rapid algorithm that can sustain the required accuracy for input parameter $y \ge 10^{-8}$ \cite{Schreier2018}.

One of the possible ways to overcome these problems is to segment the complex plane as a narrow band along $x$-axis and to use appropriate approximation for smaller $\operatorname{Im}\left[z\right] = y$ (see for example C/C++ code \cite{Johnson2012}). However, this decelerates computation as a result of additional segmentation.

Consider a complete version of the two-domain scheme \eqref{eq_3} that includes both, the real and imaginary parts, as follows
\small
\begin{equation}\label{eq_6}
w\left(z\right) =K\left(x,y\right)+iL\left(x,y\right)\approx \left\{
\begin{aligned}
&\text{interpolation}, \qquad \frac{x^2}{27^2} + \frac{y^2}{15^2} \leqslant 1 \\
&\frac{iz/\sqrt{\pi}}{z^2-1/2}, \qquad \hspace{0.65cm}\frac{x^2}{27^2} + \frac{y^2}{15^2} > 1.
\end{aligned}\right.
\end{equation}
\normalsize
The run-time test we performed shows that for the set of input numbers $\left\{x,y \in \mathbb{R}\hspace{-0.1cm}:\,\left|x\right| + y\leq 15 \right\}$ this two-domain scheme is faster in performance by a factor about $2$ as compared to that of reported in \cite{Schreier2018} (we used the MATLAB codes built on approximations \eqref{eq_4} and \eqref{eq_6}). This is possible to achieve since interpolation utilizes a simple cubic spline instead of a rational function of high order. Therefore, this fact strongly motivated us to develop further an algorithm based on a vectorized interpolation for rapid computation of the Voigt/complex error function with accuracy that meets the requirement for the HITRAN spectroscopic applications \cite{Hill2016}.

In this work we propose a new method of algorithmic implementation for rapid computation of the integrals \eqref{eq_1}, \eqref{eq_2a} and \eqref{eq_2b} by vectorized interpolation that enables us to employ just a single domain. This approach provides both, the rapid computation and accuracy that meets the requirement for the HITRAN spectroscopic applications. To the best of our knowledge, this method of computation of the Voigt/complex error function is new and has never been reported in scientific literature.

\section{Algorithmic implementation}

Computation of the Voigt function \eqref{eq_1} by interpolation with help of the lookup table was first reported in the work \cite{Drummond1985} (see also \cite{Sparks1997}). However, the lookup table involved in this method requires two dimensional interpolating grid-points that need extra memory and time to handle large-size data during computation to obtain a reasonable accuracy. In contrast, the vectorized (one dimensional) approach, where the interpolating grid points are computed dynamically at ${\mathbf{x}} = \left\{ {{x_1},{x_2},{x_3}, \ldots } \right\}$ and given fixed value $y$ \cite{Abrarov2009} (rather than picked up from the lookup table \cite{Drummond1985, Sparks1997}) is considerably advantageous in interpolation due to significantly smaller quantity of the interpolating grid-points stored in computer memory.

Let us consider two-domain scheme shown by equation \eqref{eq_6} more closely. We note that semi-major and semi-minor axises may be chosen arbitrarily depending on algorithmic implementations (see for example \cite{Poppe1990a, Poppe1990b}). Since semi-major and semi-minor axises may be flexible in magnitude, we may suggest to increase them to cover the entire set $x,y \in \mathbb{R}$ for the HITRAN applications such that we could exclude completely the rational approximation of low order in equation \eqref{eq_6} from the consideration.

Previously it was reported that the HITRAN spectroscopic database requires a domain $0 \leq \left|x\right| < 40000$ and $10^{-4} < y < {10^2}$ \cite{Wells1999}. Consequently, within the ${\text{I-st}}$ quadrant a single domain should be large enough to cover rectangular area $40000 \times 100$. Furthermore, this area should be extended by factor of four if we want to include all $4$ quadrants. However, as it has been mentioned earlier the inclusion of the all small values $y$ into consideration would be preferable for practical applications to account for low pressure and high temperature of the HITRAN gases. In our approach, the single domain represents an area $\left| x \right| \leqslant 50000$ without a boundary for the parameter $y$ since it is a scalar. Formally stating, this single domain is inclosed by a reshaped ellipse in such a way that
\[
\frac{x^2}{50000^2}+\lim_{n\to\infty}\frac{y^2}{n^2}\leqslant 1.
\]

In our work \cite{Abrarov2009} we noted that the interpolating grid-points may not be necessarily spaced equidistantly. This provides a significant advantage since we can minimize the number of the interpolating grid-points by putting them denser in more curved subintervals and sparser in more linear subintervals. Particularly, we separated interval $\left[ {-50000,50000} \right]$ along $x$-axis into subintervals as it is shown in the Table 1. Initially each subinterval contains a number of the interpolating grid-points multiple to $100$ accumulating $1600$ interpolating grid-points in total (see command lines below $<${\sffamily function MF = mainF(x,y,opt)}$>$ in the MATLAB code shown in Appendix A). However, after exclusion of repeating values the total number of the interpolating grid-points decreases from $1600$ to $197+398+4\times 198+200=1587$ (see Table 1). This small quantity of the interpolating grid-points needs a very minor amount of computer memory. Therefore, it can be easily handled by practically any modern computer. We cannot infer that such a distribution provides the smallest number of the interpolating grid-points. Perhaps, this number of the interpolating grid-points can be reduced significantly by further optimization.

\begin{table}[ht]
  \begin{center}
		  \begin{tabular}{l|r}
           \sffamily{\textbf{Subintervals}} & \sffamily{\textbf{Number of IGP}} \\
           \hline
           $\left( { - 2.5,2.5} \right)$                                          & $1\times 197=197$ \\
           $\left( { - 5.5, - 2.5} \right] \cup \left[ {2.5,5.5} \right)$         & $2\times 199=398$ \\
           $\left( { - 15, - 5.5} \right] \cup \left[ {5.5,15} \right)$           & $2\times 99=198$ \\
					 $\left( { - 100, - 15} \right] \cup \left[ {15,100} \right)$           & $2\times 99=198$ \\
					 $\left( { - 1000, - 100} \right] \cup \left[ {100,1000} \right)$       & $2\times 99=198$ \\
					 $\left( { - 10000, - 1000} \right] \cup \left[ {1000,10000} \right)$   & $2\times 99=198$ \\
					 $\left[ { - 50000, - 10000} \right] \cup \left[ {10000,50000} \right]$ & $2\times 100=200$ \\
		  \end{tabular}\vspace{0.4cm}
    \text{\sffamily{\textbf{Table 1.} Number of interpolating grid-points (IGP) in subintervals.}}
  \end{center}
\end{table}

If hypothetically there could be some input values $\left|x\right|$ greater than $50000$, then it is very easy to extend the region along $x$-axis if required. For example, it is sufficient to include only $300$ additional interpolating grid-points to extend the range for $\left|x\right|$ up to $100000$; it is not necessary to include many interpolating grid-points since the functions $K\left(x,y\right)$ and $L\left(x,y\right)$ become nearly linear as $\left|x\right|$ increases.

For proper interpolation the accuracy of the computation should be two orders of the magnitude better than $10^{-6}$. Generally, any highly accurate algorithm can be used to compute $1587$ interpolating grid-points. We applied the highly accurate MATLAB function {\sffamily{fadsamp.m}} shown in our recent paper \cite{Abrarov2018a} as a subroutine for this purpose (see a brief description of this method in Appendix B). Alternatively, a MATLAB function file that is highly accurate and suitable for computation of the interpolating grid-points can be found in \cite{Abrarov2011} (these two codes can also be downloaded from the Matlab Central websites \cite{MatlabCentral66752} and \cite{MatlabCentral47801}, respectively). Highly accurate C/C++ implementation by Johnson \cite{Johnson2012} with MEX plugins for MATLAB \cite{MatlabCentral38787} can also be used for a subroutine that can be invoked from the MATLAB environment to compute these $1587$ interpolating grid-points.

It should be noted that this methodology can also be generalized further to other functions (including spectral line profiles). For example, our preliminary results demonstrate that this methodology is also applicable for rapid and quite accurate computation of the Spectrally Integrated Voigt Function (SIVF) that enables us to perform line-by-line atmospheric modeling at reduced spectral resolution \cite{Quine2013}.

\section{Error analysis and run-time test}

In order to perform error analysis define the relative errors for the real
$$
{\Delta _{\operatorname{Re} }} = \left| {\frac{{{K_{{\text{ref}}{\text{.}}}}\left( {x,y} \right) - K\left( {x,y} \right)}}{{K_{{\text{ref}}{\text{.}}}{{\left( {x,y} \right)}}}}} \right|
$$
and imaginary parts
$$
{\Delta _{\operatorname{Im} }} = \left| {\frac{{{L_{{\text{ref}}{\text{.}}}}\left( {x,y} \right) - L\left( {x,y} \right)}}{{L_{{\text{ref}}{\text{.}}}{{\left( {x,y} \right)}}}}} \right|
$$
of the complex error function \eqref{eq_1}, where ${K_{{\text{ref}}{\text{.}}}}\left( {x,y} \right)$ and ${L_{{\text{ref}}{\text{.}}}}\left( {x,y} \right)$ are the highly accurate reference values that can be readily obtained by using the CAS.

\begin{figure}[ht]
\begin{center}
\includegraphics[width=32pc]{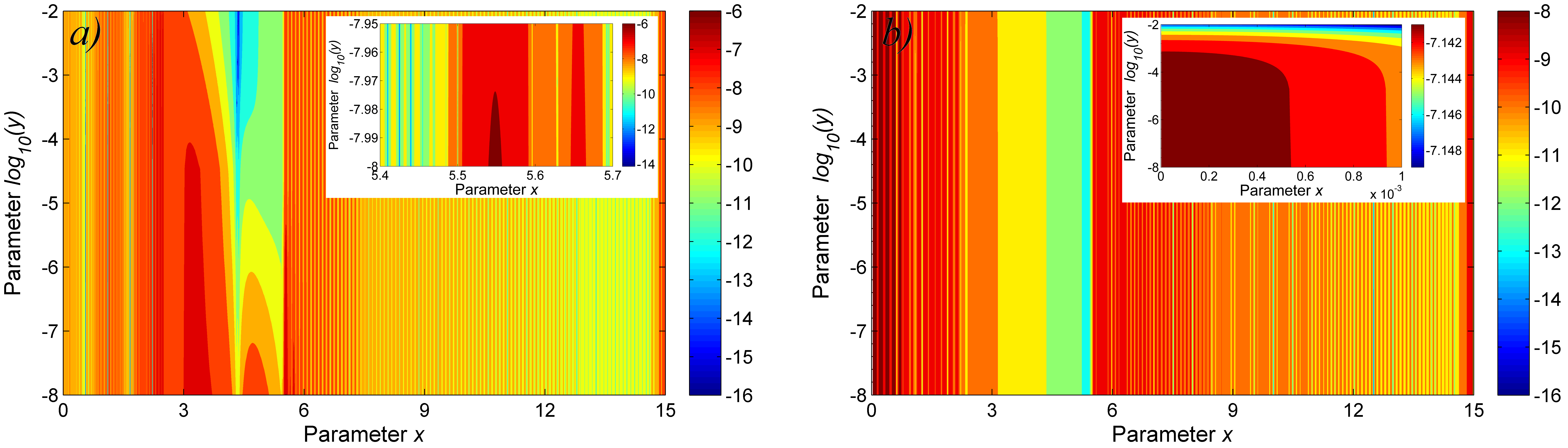}\hspace{2pc}%
\begin{minipage}[b]{28pc}
\vspace{0.3cm}
{\sffamily {\bf{Fig. 1.}} Logarithms of relative error for the real a) and imaginary b) parts of the complex error function over the area $x\in\left[0,15\right]\cap y\in\left[10^{-8},10^{-2}\right]$ computed by vectorized interpolation. Insets show the subareas with worst accuracies for the real a) and imaginary b) parts.}
\end{minipage}
\end{center}
\end{figure}

Figures 1a and 1b show the relative errors for the real and imaginary parts over the area $x \in \left[ {0,15} \right]$ and $y\in\left[10^{-8},10^{-2}\right]$ computed by vectorized interpolation. The largest relative errors over this area for the real and imaginary parts are $1.0589\times 10^{ - 6}$ and $7.236 \times {10^{ - 8}}$, respectively. Thus, our algorithm satisfies the accuracy criterion $10^{-6}$ for the required range $y\geq 10^{-8}$. In comparison, the computational test we perform reveals that approximation \eqref{eq_4} sustains the required accuracy only at $y \ge 10^{-6}$.

Figures 2a and 2b show the relative errors for the real and imaginary parts over the area $x \in \left[ {0,15} \right]$ and $y \in \left[ {{{10}^{ - 2}},15} \right]$. The largest relative errors over this area for the real and imaginary parts are $2.7766\times {10^{ - 7}}$ and $7.0619\times {10^{ - 8}}$, respectively.

Thus, from Figs. 1a, 1b and 2a, 2b we can conclude that our algorithm satisfies accuracy requirement for the HITRAN-based applications and effectively resolves the problem that is inherent to the single quotient shown in equations  \eqref{eq_4} and \eqref{eq_5} as well as many other approximations \cite{Armstrong1967, Amamou2013}.

In order to obtain most objective results for the run-time test we used an independently written code for the Voigt/complex error function. Specifically, the run-time test has been performed by comparing our MATLAB code shown in Appendix A with C/C++ implementation that was developed by Johnson \cite{Johnson2012, MatlabCentral38787}. This implementation is known to be the fastest C/C++ program. It represents a modified Algorithm $680$ \cite{Poppe1990a, Poppe1990b} with inclusion of the Salzer's approximations \cite{Abrarov2018b} (see also \cite{Salzer1951}). As the computation complexity prevails at smaller values of the parameter $y$, we imply that it is close to zero, say $y=10^{-5}$. The detailed description of how to run the programs and perform the time execution test is shown in Appendix C.

The run-time test shows that for $10$ million input values within the most important area such that $\left\{x,y\in\mathbb{R}\hspace{-0.1cm}:\left|x+iy\right|\leq 15\right\}$, the MATLAB code is almost $8$ times faster than the C/C++ implementation while for the entire domain that is required for the HITRAN spectroscopic applications the MATLAB code is faster than the C/C++ implementation by a factor about $3$.

\begin{figure}[ht]
\begin{center}
\includegraphics[width=32pc]{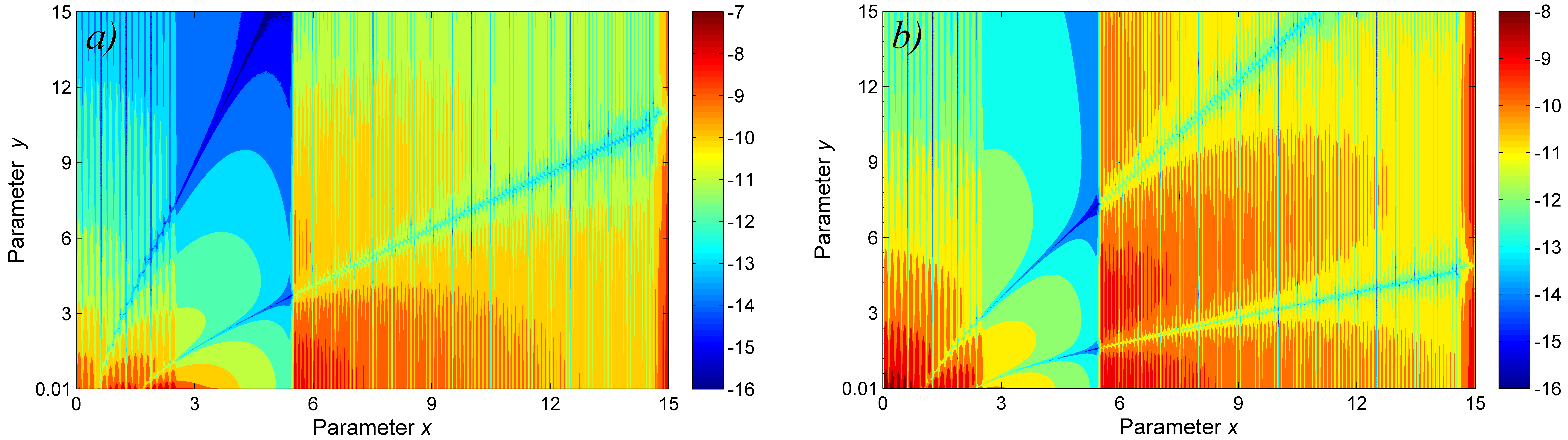}\hspace{2pc}%
\begin{minipage}[b]{28pc}
\vspace{0.3cm}
{\sffamily {\bf{Fig. 2.}} Logarithms of relative error for the real a) and imaginary b) parts of the complex error function over the area $x\in\left[0,15\right]\cap y\in\left[10^{-2},15\right]$ computed by vectorized interpolation.}
\end{minipage}
\end{center}
\end{figure}

The MATLAB is one of the fastest scientific languages in computation. However, it is generally slower than C/C++. The more rapid computation has been achieved because of two main reasons. The function {\sffamily{Faddeeva.cc}} is unnecessarily complicated as it utilizes a large number of domains that due to multiple logical operations and sorting of the input numbers $x$ and $y$ decelerate computation. Furthermore, our approach {\it{de facto}} performs $1587$ actual computations only as the remaining is just an interpolation. By choosing different options for interpolation, we found experimentally that the MATLAB built-in method {\sffamily'spline'} provides the best performance. All these results can be readily confirmed by running the codes provided in Appendices A and C.

\section{Conclusion}

We propose a new single-domain vectorized interpolation method for rapid computation of the Voigt/complex error function \eqref{eq_1} that enables us to achieve required accuracy for the HITRAN-based spectroscopic applications. The computational test we performed reveals that within intervals $x\in\left[0,15\right]\cap y\in\left[10^{-8},15\right]$ and $x\in\left[0,50000\right] \cap y \geq 10^{-8}$ our algorithmic implementation is faster in computation by factors of about $8$ and $3$, respectively, as compared to the fastest known C/C++ code for the Voigt/complex error function.

\section*{Acknowledgments}

This work is supported by National Research Council Canada, Thoth Technology Inc., York University, Epic College of Technology and Epic Climate Green (ECG) Inc. The authors wish to thank principal developer of the MODTRAN Dr. Alexander Berk for constructive discussions.

\newpage
\bigskip
\section*{Appendix A}

\footnotesize
\begin{verbatim}

function vecFF = vecfadf(x,y,opt)

% This function file computes the Voigt/complex error function, also known
% as the Faddeeva function, providing rapid computation at required
% accuracy for the HITRAN-based radiative transfer application.
%
% SYNOPSIS:
%          x   - row or column vector
%          y   - scalar
%          opt - option for the real and imaginary parts
%
% The code is written by Sanjar M. Abrarov, Brendan M. Quine, Rehan
% Siddiqui and Rajinder K. Jagpal, York University, Canada, May 2019.

bound = 5*1e4; % default bound to cover the HITRAN domain
num = 1e2; % common number for grid-points in interpolation

if max(size(y)) ~= 1
    disp('Parameter y must be a scalar!')
    return
elseif ~isvector(x)
    disp('Parameter x is not a vector!')
end

if max(abs(x)) > bound || abs(y) < 1e-8
    disp('x or y is beyond HITRAN range! Computation is terminated.')
    return
end

if nargin == 2
    opt = 3;
    disp('Default value opt = 3 is assigned.')
end

if opt ~= 1 && opt ~= 2 && opt ~= 3
    disp(['Wrong parameter opt = ',num2str(opt),'! Use either 1, 2 or 3.'])
    return
end

if y >= 0
    vecFF = mainF(x,y,opt); % upper half-plane
else
    vecFF = mainF(x,-y,opt); % lower half-plane
    vecFF = conj(2*exp(-(x + 1i*y).^2) - vecFF);
end

    function MF = mainF(x,y,opt)

        % Forming non-equidistantly spaced interpolating grid-points (IGP)
        IGP = linspace(0,2.5,num);
        IGP = [IGP,linspace(2.5,5.5,100 + num)]; % 100 more grid-points
        IGP = [IGP,linspace(5.5,15,num)];
        IGP = [IGP,linspace(15,100,num)];
        IGP = [IGP,linspace(100,1000,num)];
        IGP = [IGP,linspace(1000,10000,num)];
        IGP = [IGP,linspace(10000,bound,num)];
        IGP = [-(flip(IGP)),IGP];
        IGP = unique(IGP); % exclude repeated values

        MF = interp1(IGP,fadsamp(IGP + 1i*y),x,'spline'); % call ...
        % external function <fadfsamp.m>. This MATAB function file is ...
        % provided in paper [Abrarov, Quine & Jagpal, Appl. Num. Math., ...
        % 129 (2018) 181-191]. 
        % URL: https://doi.org/10.1016/j.apnum.2018.03.009
                
        switch opt
            case 1
                MF = real(MF);
            case 2
                MF = imag(MF);
        end
    end
end

\end{verbatim}
\normalsize

\newpage
\bigskip
\section*{Appendix B}

In our publications \cite{Quine2013, Abrarov2015a} we have introduced the following product-to-sum identity
\[
\prod_{m=1}^{k}\cos\left(\frac{t}{2^m}\right)=\frac{1}{2^{k-1}}\sum_{m=1}^{2^{k-1}}\cos\left(\frac{m-1/2}{2^{k-1}}t\right), \qquad \forall k\geq 1
\]
and since \cite{Gearhart1990, Weisstein3}
\[
\text{sinc}\left(t\right)=\prod_{m=1}^\infty\cos\left(\frac{t}{2^m}\right)
\]
from this product-to-sum identity it immediately follows that the sinc function can be expanded as a sum of cosines 
\[
\text{sinc}\left(t\right)=\lim_{k\to\infty}\frac{1}{k}\sum_{m=1}^{k}\cos\left(\frac{m-1/2}{k}t\right)
\]
or
\[\tag{A.1}\label{eq_A.1}
\text{sinc}\left(t\right)\approx\frac{1}{k}\sum_{m=1}^{k}\cos\left(\frac{m-1/2}{k}t\right), \qquad k >> 1.
\]

Using a new method of sampling based on incomplete cosine expansion of the sinc function \eqref{eq_A.1} we can obtain \cite{Abrarov2015a} (see also \cite{Abrarov2015b} and cited literature in context therein)
\[\tag{A.2}\label{eq_A.2}
\begin{aligned}
 w\left( z \right) \approx & \,\Omega \left( {z + i\varsigma /2} \right) \\ 
 \Rightarrow & \,\Omega \left( z \right) \triangleq \sum\limits_{m = 1}^M {\frac{{{A_m} + {B_m}z}}{{C_m^2 - {z^2}}}}, 
\end{aligned}
\]
where $N = 23$, $M=23$, $h=0.25$, $\varsigma =2.75$ and the expansion coefficients
$$
{A_m} = \frac{{\sqrt \pi  \left( {m - 1/2} \right)}}{{2{M^2}h}}\sum\limits_{n =  - N}^N {{e^{{\varsigma ^2}/4 - {n^2}{h^2}}}\sin \left( {\frac{{\pi \left( {m - 1/2} \right)\left( {nh + \varsigma /2} \right)}}{{Mh}}} \right)},
$$
$$
{B_m} =  - \frac{i}{{M\sqrt \pi  }}\sum\limits_{n =  - N}^N {{e^{{\varsigma ^2}/4 - {n^2}{h^2}}}\cos \left( {\frac{{\pi \left( {m - 1/2} \right)\left( {nh + \varsigma /2} \right)}}{{Mh}}} \right)} 
$$
and
$$
{C_m} = \frac{{\pi \left( {m - 1/2} \right)}}{{2Mh}}.
$$

Similar to the Huml\'{i}\v{c}ek's approximation \eqref{eq_5}, the series expansion \eqref{eq_A.2} also deteriorates in accuracy with decreasing $y$. We have shown how to overcome this problem by transformation of equation \eqref{eq_A.2} into following form
\[\tag{A.3}\label{eq_A.3}
w\left( z \right) \approx e^{-z^2} + z\sum\limits_{m = 1}^{M + 2} {\frac{{{\kappa _m} - {\lambda _m}{z^2}}}{{{\mu _m} - {\nu _m}{z^2} + {z^4}}}},
\]
where
$$
{\kappa _m} = {B_m}\left[C_m^2 - \left(\frac{\varsigma ^2}{2}\right)^2 \right] + i{A_m}\varsigma = {B_m}\left[ {{{\left( {\frac{{\pi \left( {m - 1/2} \right)}}{{2Mh}}} \right)}^2} - {{\left( {\frac{\varsigma }{2}} \right)}^2}} \right] + i{A_m}\varsigma,
$$
$$
{\lambda_m} = {B_m},
$$
$$
{\mu _m} = C_m^4 + \frac{C_m^2{\varsigma ^2}}{2} + \frac{\varsigma ^4}{16} = {\left[ {{{\left( {\frac{{\pi \left( {m - 1/2} \right)}}{{2Mh}}} \right)}^2} + {{\left( {\frac{\varsigma }{2}} \right)}^2}} \right]^2}
$$
and
$$
{\nu _m} = 2C_m^2 - \frac{\varsigma ^2}{2} = 2{\left( {\frac{{\pi \left( {m - 1/2} \right)}}{{2Mh}}} \right)^2} - {\frac{\varsigma }{2}^2}.
$$

The third equation that is used in the function file {\sffamily{fadsamp.m}} is the Laplace continued fraction \cite{Abramowitz1972, Gautschi1970, Poppe1990a} given by
\[\tag{A.4}\label{eq_A.4}
w\left( z \right) \approx \frac{{\left( {i/\sqrt \pi  } \right)}}{{z - \frac{{1/2}}{{z - \frac{1}{{z - \frac{{3/2}}{{z - \frac{2}{{z - \frac{{5/2}}{{z - \frac{3}{{z - \frac{{7/2}}{{z - \frac{4}{{z - \frac{{9/2}}{{z - \frac{5}{{z - \frac{{11/2}}{z}}}}}}}}}}}}}}}}}}}}}}},
\]

Equations \eqref{eq_A.2} and \eqref{eq_A.3} cover the domains
$$
\left\{x,y\in \mathbb{R}\hspace{-0.1cm}:\left|x+iy \right|\leq 8 \right\}\setminus \left\{y<0.05\left|x\right|\right\}$$
and
$$
\left\{x,y\in \mathbb{R}\hspace{-0.1cm}:\left|x+iy \right|\leq 8\right\}\cap \left\{y<0.05\left|x\right|\right\},
$$
respectively, while equation \eqref{eq_A.3} covers the domain $\left\{x,y\in \mathbb{R}\hspace{-0.1cm}:\left|x+iy \right| > 8 \right\}$ (see Fig. 1 in \cite{Abrarov2018a}). This three-domain scheme excludes all poles in computation and provides largest relative error $\sim 10^{-13}$ only.

\newpage
\section*{Appendix C}

The C/C++ function file {\sffamily{Faddeeva.cc}} and its header {\sffamily{Faddeeva.hh}} can be downloaded from the website \cite{Johnson2012}. To run the program the following lines can be added to the end of the {\sffamily{Faddeeva.cc}} function file

\vspace{-0.5cm}
\footnotesize
\begin{verbatim}

/**
These command lines can be added at the end of the source
file 'Faddeeva.cc'.
*/

#ifdef __cplusplus
# include <cstdio>
# include <cstdlib>
# include <iostream>
# include <iomanip>

#else
# include <stdio.h>
#endif

/**
The function fRand(minNum, maxNum) generates a random number
within the range from 'minNum' to 'maxNum'.
*/

double fRand(double minNum, double maxNum){

    double fMin = minNum, fMax = maxNum;
    double f = (double)rand() / RAND_MAX;

return fMin + f * (fMax - fMin);
};

int main(void){

double(epsVal) = 1E-6; // epsilon value for accuracy

    for(int k = 0; k < 1e7; k++){ // execute the computation ...
    // 10 million times with random numbers 0 < x < 15 and positive y << 1.

        FADDEEVA(w)(C(fRand(0,15),1E-5),epsVal);
    };
return 0;
};

\end{verbatim}
\normalsize

As we can see from the body of the function \textless{\sffamily int main(void)}\textgreater, this program computes $10$ million random digits of $x$ at fixed value of $y = 10^{-5}$. The epsilon value that determines accuracy of computation is taken to be $10^{-6}$. The execution time on a typical laptop computer (we used Intel(R) Core(TM) i3-7020U CPU @ 2.30GHz, 8GB RAM) takes about $14$ seconds.

It is interesting to note that with \textless{\sffamily double(espVal) = 1E-12;}\textgreater\, the executions time increases to $17$ seconds. As a subsequent step, at \textless{\sffamily double(espVal) = 0;}\textgreater\, we should expect the highest accuracy and, therefore, the longest execution time. Instead, however, the program becomes significantly faster and it takes about $8.5$ seconds to compute all $10$ million digits. This unexpected behavior of the C/C++ implementation suggests that the program, most likely, does not guarantee the prescribed accuracy for all input numbers $x$ and $y$ at \textless{\sffamily double(espVal) = 0;}\textgreater. Therefore, we did not take this option as a reference.

The following is the MATLAB code that also computes $10$ million random numbers $x \in \left[0,15\right]$ at fixed $y=10^{-5}$. The execution time is about $1.8$ seconds only. As we can see, the MATLAB code based on interpolation is almost $8$ times faster than the C/C++ implementation.

\vspace{-0.5cm}
\footnotesize
\begin{verbatim}

% ******************************************************************
x = 15*rand(1e7,1); % this generates 10 million random numbers ...
                    % in the interval 0 < x < 15
y = 1e-5; % y is a scalar
tic;vecfadf(x,y,3);toc % this shows the run-time
% ******************************************************************

\end{verbatim}
\normalsize

\vspace{-0.5cm}
In order to compare the execution times for the interval $x \in \left[0,50000\right]$, we can simply replace the command lines above in the C/C++ and MATLAB codes as
\vspace{-0.2cm}
\footnotesize
\begin{verbatim}
FADDEEVA(w)(C(fRand(0,50000),1E-5),epsVal);
\end{verbatim}
\normalsize
and
\vspace{-0.2cm}
\footnotesize
\begin{verbatim}
x = 50000*rand(1e7,1); % this generates 10 million random numbers ...
                       % in the interval 0 < x < 50000
\end{verbatim}
\normalsize
respectively. The execution times for the C/C++ and MATLAB implementations for this case become about $3$ and about $1$ seconds, respectively. Therefore, the MATLAB code is nearly $3$ times faster.

Execution times can be decreased by more than an order of the magnitude on a more powerful computer of the latest generation. However, we should anticipate that these ratios $8$ and $3$ will remain same.


\end{document}